\documentclass{ip-journal}
\usepackage{amsmath}
\usepackage{amsthm}
\usepackage{amsfonts}
\usepackage{amssymb}
\usepackage{mathrsfs}
\usepackage{tikz}
\usepackage{environ}
\usepackage{elocalloc}
\usepackage{url}
\usepackage{caption}
\usepackage{CJKutf8}
\usepackage{mathtools}
\usepackage{dcpic}
\usepackage{tikz-cd}
\usepackage{stmaryrd}
\allowdisplaybreaks
\usetikzlibrary{positioning}
\usetikzlibrary{arrows,chains,positioning,scopes,quotes}
\usetikzlibrary{decorations.markings}

\title{Stability Inequalities for Lawson Cones}
\author{Zhenhua Liu}
\date{}
\dedicatory{Dedicated to Xunjing Wei}

\begin{document}
	\maketitle
	\newcommand{\ai}{\alpha}
	\newcommand{\be}{\beta}
	\newcommand{\Ga}{\Gamma}
	\newcommand{\ga}{\gamma}	
	\newcommand{\de}{\delta}
	\newcommand{\De}{\Delta}
	\newcommand{\e}{\epsilon}
	\newcommand{\lam}{\lambda}
	\newcommand{\La}{\Lamda}
	\newcommand{\om}{\omega}
	\newcommand{\Om}{\Omega}
	\newcommand{\si}{\sigma}
	\newcommand{\vp}{\varphi}
	\newcommand{\rh}{\rho}
	\newcommand{\ta}{\theta}
	\newcommand{\Ta}{\Theta}
	\newcommand{\W}{\mathcal{O}}
	\newcommand{\ps}{\psi}
	
	\newcommand{\mf}[1]{\mathfrak{#1}}
	\newcommand{\ms}[1]{\mathscr{#1}}
	\newcommand{\mb}[1]{\mathbb{#1}}
	\newcommand{\cd}{\cdots}
	
	\newcommand{\s}{\subset}
	\newcommand{\es}{\varnothing}
	\newcommand{\cp}{^\complement}
	\newcommand{\bu}{\bigcup}
	\newcommand{\ba}{\bigcap}
	\newcommand{\co}{^\circ}
	\newcommand{\ito}{\uparrow}
	\newcommand{\dto}{\downarrow}

	\newcommand{\du}{^\ast}
	\newcommand{\pf}{_\ast}
	\newcommand{\m}{^{-1}}
	\newcommand{\ts}{\otimes}
	\newcommand{\ip}{\cdot}

	\newcommand{\w}{\wedge}
	\newcommand{\X}{\mathfrak{X}}
	\newcommand{\pd}{\partial}
	\newcommand{\dx}{\dot{x}}
	\newcommand{\dr}{\dot{r}}
	\newcommand{\dy}{\dot{y}}
	\newcommand{\dth}{\dot{theta}}
	\newcommand{\pa}[2]{\frac{\pd #1}{\pd #2}}
	\newcommand{\na}{\nabla}
	\newcommand{\dt}[1]{\frac{d#1}{d t}\bigg|_{ t=0}}
	\newcommand{\ld}{\mathcal{L}}

	\newcommand{\N}{\mathbb{N}}
	\newcommand{\R}{\mathbb{R}}
	\newcommand{\Z}{\mathbb{Z}}
	\newcommand{\Q}{\mathbb{Q}}
	\newcommand{\C}{\mathbb{C}}
	\newcommand{\bh}{\mathbb{H}}
	\newcommand{\hi}{\mathcal{H}}
	\newcommand{\lb}{\mathcal{L}}

	\newcommand{\II}{\textnormal{II}}
	\newcommand{\di}{\operatorname{div}}
	\newcommand{\dist}{\textnormal{dist}}

	\newcommand{\wh}{\Rightarrow}
	\newcommand{\eq}{\Leftrightarrow}
	
	\newcommand{\ho}{\textnormal{Hom}}
	\newcommand{\ds}{\displaystyle}

	\theoremstyle{plain}
	\newtheorem{thm}{Theorem}[subsection]
	\renewcommand{\thethm}{\arabic{thm}}
		\newtheorem{resu}{Result}

	\newtheorem{lem}[thm]{Lemma}
	\newtheorem{prop}[thm]{Proposition}
	\newtheorem*{cor}{Corollary}
	\newtheorem*{pro}{Proposition}
	
	\theoremstyle{definition}
	\newtheorem{defn}{Definition}[section]
	\newtheorem{conj}{Conjecture}[section]
	\newtheorem{exmp}{Example}[section]
	
	\theoremstyle{remark}
	\newtheorem*{rem}{Remark}
	\newtheorem*{note}{Note}
	
	\begin{abstract}
In paper \cite{phi}, G. De. Philippis and F. Maggi proved global quadratic stability inequalities and derived explicit lower bounds for the first eigenvalues of the stability operators for all area-minimizing Lawson cones $M_{kh}$, except for those with
		$$(k,h),(h,k)\in S=\{(3,5),(2,7),(2,8),(2,9),(2,10),(2,11)\}.$$
		We proved the corresponding inequalities and lower bounds for these Lawson cones $M_{kh}$ with $(k,h),(h,k)\in S$ by using different sub-calibrations from theirs, thus extending their results to all area-minimizing Lawson cones.
	\end{abstract}
	\section{Introduction}
Suppose $h,k\ge 2$ are positive integers. The Lawson cone $M_{kh}$ is the level set
	\begin{align*}
	M_{kh}=\left\{z=(x,y)\in \R^k\times \R^h:\frac{|x|}{\sqrt{k-1}}=\frac{|y|}{\sqrt{h-1}}\right\}.
	\end{align*}
It is known to be area-minimizing  (see \cite{c1}, \cite{c2}, \cite{c3}, and \cite{c4}) provided
	\begin{align}\label{lc}
	h+k\ge 9,\text{ or }(h,k)=(3,5),(4,4),(5,3).
	\end{align}	
In their paper \cite{phi}, G. De. Philippis and F. Maggi proved  global quadratic stability inequalities and derived explicit lower bounds for the first eigenvalues of the stability operators for all area-minimizing Lawson cones $M_{kh}$, except for 
\begin{align*}
(h,k),(k,h)\in S=\{&(3,5),(2,7),(2,8),(2,9),(2,10),(2,11)\}.
\end{align*}
They achieved this by exploiting sub-calibrations for Lawson cones. Unfortunately, the sub-calibrations that they used did not work for the cones $M_{kh}$ with $(h,k),(k,h)\in S.$
	 Our main results, Theorem 1 and Theorem 2 in Section 1.1, extend these inequalities to the cones $M_{kh}$ with $(h,k),(k,h)\in S.$ We achieve this by carefully choosing sub-calibrations for these Lawson cones in Lemma 2 of Section 2.1. However, our sub-calibrations do not work for other cases in general.
	 	 
	 We first review their results and explain their methods, which we mostly follow. Consider a variation with compact support of the Lawson cone $M_{kh}.$ Suppose the variation can be realized as the boundary of a set $F$ of finite perimeter. Roughly speaking, their first result controls the volume bounded between the Lawson cone and the variation $\pd F$ by the difference between the area of the variation $\pd F$ and that of the cone $M_{kh}$ up to scaling. Their second result provides lower bounds for the first eigenvalues of the stability operators. For a great discussion of the significance of these results, please refer to Section 1 of \cite{phi}.
	 
The Lawson cone $M_{kh}$ can be realized as the boundary $\pd K_{kh}$ of the region	
\begin{align*}
K_{kh}=&\left\{(x,y)\in \R^k\times \R^h:\frac{|x|}{\sqrt{k-1}}<\frac{|y|}{\sqrt{h-1}}\right\}.
\end{align*}	
Let $\lb^m$ denote the $m$-dimensional Lebesgue measure, $\om_n$ denote the volume of unit $n$-ball, and $P(A;B)$ denote the perimeter of $A$ in $B.$ Their results are as follows.
\begin{resu}(Theorem 5 in \cite{phi})
	If $R>0,m=h+k,(h,k)\not\in S$ satisfy all the conditions in (1),  then
	\begin{align*}
	\left(\frac{\lb^m(K_{kh}\De F)}{R^m}\right)^2\le C\frac{P(F;H_R)-P(K_{kh}; H_R)}{R^{m-1}},
	\end{align*}whenever $F$ is a set of locally finite perimeter with symmetric difference $K_{kh}\De F\s \s H_R= B_R^k\times B_R^h.$ 
	Possible values of $C$ are
	\begin{align*}
	C=&\frac{2^{12}\sqrt{\om_k\om_h}}{(k-1)^{1/8}}\sqrt{\frac{hk}{m-1}}(\frac{h-1}{k-1})^{3/2},\text{ if }2\le k\le h, (k,h)\not=(4,4),\\&\text{Interchange }k,h\text{ if }2\le h\le k.\\
C=&128\om_4,\text{ if }(k,h)=(4,4).
	\end{align*}
\end{resu}

\begin{resu}(Theorem 2 in \cite{phi})
	If $R,m,h,k$ are as in Result 1, and
	\begin{align*}
	\lam_{k,h}(R)=\inf\bigg\{\int_{M_{kh}}^{}|\na^{M_{kh}}\vp|^2-|\textnormal{\II}_{{M_{kh}}}|^2\vp^2 d\mathcal{H}^{m-1}:\int_{M_{kh}}\vp^2=1,\textnormal{spt}\vp\s\s B^m_R\bigg\},
	\end{align*}then	
	\begin{align*}
	\lam_{k,h}(R)\ge \frac{c_{k,h}}{R^2}.
	\end{align*}
	Possible values of $c_{k,h}$ are 
		\begin{align*}
		c_{k,h}=&\frac{1}{2^9}\bigg(\frac{k-1}{h-1}\bigg)^{9/4}\frac{(m-2)^{1/2}}{(h-1)^{1/4}},\text{ if }2\le k\le h, (k,h)\not=(4,4).\\&\text{Interchange }k,h\text{ if }2\le h\le k.\\
		c_{k,h}=&\frac{\sqrt{2}}{16},\text{ if }(k,h)=(4,4).
		\end{align*}
\end{resu}
As illustrated in Figure 1, their method is based on sub-calibrating the Lawson cones with a unit-length vector field $g.$ In other words, the vector field $g$ restricts to the unit normal on $M_{kh}$, and the divergence $\di g$ does not change sign in $K_{kh}$ and $K_{kh}\cp$, respectively. 
\begin{figure}[h]
\centering
\includegraphics[width=1\linewidth]{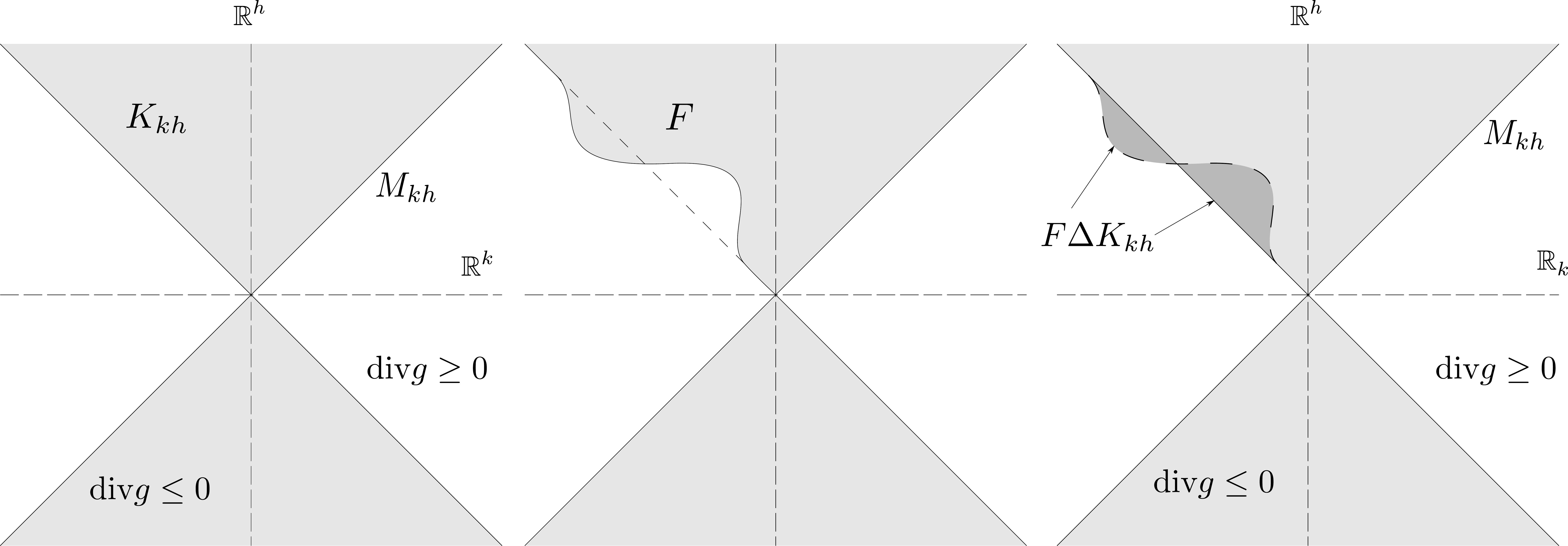}
\caption{A sub-calibration $g$ of the Lawson cone $M_{kh}$ and a variation.}
\label{fig:drawing}
\end{figure}

After cleverly choosing $g$, they proved that 
\begin{align}
\di  g(z)\ge c_{k,h}\frac{\dist(z,M_{kh})}{|z|^2},
\end{align}where $\dist$ is the Euclidean distance. Then they exploit inequality (2) to deduce the desired results. For a beautiful discussion of sub-calibrations (also called quantitative calibrations), please refer to their paper \cite{phi}. 

Unfortunately, the sub-calibrations they used did not work for $(h,k),(k,h)\in S$. The main results of this paper extend their stability inequalities to include those $(k,h).$ We achieve this by using sub-calibrations inspired by \cite{c4}.

\subsection{Stability Inequalities Extended to $(h,k),(k,h)\in S$}
	 \begin{thm}
	 	If $R>0,m=h+k,(h,k),(k,h)\in S,$ then
	 	\begin{align*}
	 	\left(\frac{\lb^m(K_{kh}\De F)}{R^m}\right)^2\le C\frac{P(F;H_R)-P(K_{kh}; H_R)}{R^{m-1}},
	 	\end{align*}whenever $F$ is a set of locally finite perimeter with $K_{kh}\De F\s \s B_R^k\times B_R^h.$ 
	 	A possible value of $C$ is $7^2\times 12^2\times 10^{20}.$
	 \end{thm}

	 \begin{thm}
	 	If $R,m,h,k$ are as in Theorem 1, and
	 	\begin{align*}
	 	\lam_{k,h}(R)=\inf\bigg\{\int_{M_{kh}}^{}|\na^{M_{kh}}\vp|^2-|\textnormal{\II}_{{M_{kh}}}|^2\vp^2 d\mathcal{H}^{m-1}:\int_{M_{kh}}\vp^2=1,\textnormal{spt}\vp\s\s B^m_R\bigg\},
	 	\end{align*}then	
	 	\begin{align*}
	 	\lam_{k,h}(R)\ge \frac{c_{k,h}}{R^2},
	 	\end{align*}
	 	Possible values of $c_{k,h}$ are 
	 	\begin{align*}
	 	c_{3,5}=c_{5,3}&=\frac{\sqrt{3}}{21^3},\\c_{k,2}=c_{2,h}&=\frac{\sqrt{11}}{11^6},
	 	\end{align*} for $k,h=7,8,9,10,11.$
	 \end{thm}
	
	\section{Proof of the Theorems}
	We now prove, in order, Theorem 2 and Theorem 1. By the symmetry of Lawson cones, it suffices to prove the cases with $(h,k)\in S.$ The following lemma is the basic tool to extract information from the sub-calibrations $g.$
	\begin{lem}
		If $m\ge 2,$ $E$ is of locally finite perimeter in $\R^m,$ and $g\in W^{1,1}_{\textnormal{loc}}(\R^m,\R^m),$ 
		\begin{align*}
		|g|&\le 1\text{ on }\R^m,\\
		\di g&\ge 0,\text{ a.e. on }E^c,\\
		\di g&\le 0,\text{ a.e. on }E,\\
		g&=\nu_E,\text{ }\hi^{m-1}-\text{a.e. on }\pd_{1/2}E,
		\end{align*}then $E$ is a local minimizer of the perimeter in $\R^m,$ with 	\begin{align}\label{va}
		P(F;A)-P(E;A)=\int_{E\De F}|\di g|+\int_{A\cap \pd_{1/2}F}1-(g\ip \nu_F)d\hi^{m-1}.
		\end{align}
	\end{lem}
Here $\hi^{m-1}$ is the $m-1$-dimensional Hausdorff measure, $\nu_E$ is the out-pointing unit normal. If $|E|$ denote the $\lb^m$-volume of a set $E$,  then
\begin{align*}
\pd_{1/2}E=\{x\in\R^m:\lim_{r\to 0^+}\frac{|E\cap B(x,r)|}{\om_nr^n}=\frac{1}{2}\},
\end{align*}is defined as the set of points of density $1/2$ in $E.$ For proof of Lemma 1 and details about $\pd_{1/2}E$, please refer to the proof of Proposition 4.1 in \cite{phi} and the relevant discussions on page 416 in \cite{phi}.	Roughly speaking, Lemma 1 can be proved by breaking down the integration definition of perimeter and then using the divergence theorem.
	
	The left hand-side of (\ref{va}) can be seen as variation of area, so it can provide information for second variation by Taylor expansion and choosing suitable variation $F$. The key to using this information is to find vector fields $g$ that satisfy inequality (1) in Section 1.
	\subsection{Sub-calibrations for $M_{kh}$ with $(h,k)\in S$}
	\begin{lem}
	For $E=K_{kh},$ the vector field $$g=\frac{\na f}{|\na f|}$$ satisfies all the hypothesis in Lemma 1. The function $f$ for $(h,k)=(3,5)$ is 
	\begin{align*}
	f(x,y)=\begin{cases}\ds
	\frac{(h-1)|x|^2-(k-1)|y|^2}{4}((h-1)|x|)^{3/2},&\text{if }z\in K_{kh},\\
	\ds\frac{(h-1)|x|^2-(k-1)|y|^2}{4}((k-1)|y|)^{3/2},&\text{if }z\in K_{kh}\cp,
	\end{cases}
	\end{align*}
	and the functions $f$ for $(h,k)=(2,k)$ with $k=7,8,9,10,11$ are
	\begin{align*}
	f(x,y)=\begin{cases}
	\ds\frac{(h-1)|x|^2-(k-1)|y|^2}{4}((h-1)|x|)^{3},&\text{if }z\in K_{kh},\\
	\ds\frac{(h-1)|x|^2-(k-1)|y|^2}{4}((k-1)|y|)^{2},&\text{if }z\in K_{kh}\cp.
	\end{cases}
	\end{align*}Moreover, $g$ also satisfy 
	\begin{align*}
	|\di g|\ge\frac{c_{k,h}}{|z|^2}\dist(z,M_{kh}),
		\end{align*} with values of $c_{k,h}$ the same as in Theorem 2.
	\end{lem}
The proof of Lemma 2 is left to Section 3. The sub-calibrations we choose work well for $(h,k)\in S,$ but do not work for some other Lawson cones. In some sense, these are specifically chosen to cover the cases $(h,k)\in S.$
\subsection{Proof of Theorem 2}
	By Lemma 1, we have
	\begin{align*}
	P(F;H_R)-P(K_{kh};H_R)\ge& \int_{K_{kh}\De f}|\di g|\\\ge& c_{k,h}\int_{K_{kh}\De F}\frac{\dist(z,M_{kh})}{|z|^2}dz\\\ge& \frac{c_{k,h}}{R^2}\int_{K_{kh}\De F}\dist(z,M_{kh})dz.
	\end{align*}
	Now, suppose $\vp\in C^1(M_{kh}),$ with $0\not\in\text{spt}\vp\s\s B_R^m.$ For $t_0>0$ small enough, there exists an open set $F\s \R^m$ with $\pd F-\{0\}$ a $C^1$ hypersurface and $K_{kh}\De F\s\s H_R,$ such that
	\begin{align*}
	\pd F-\{0\}=\{z+t\vp(z)\nu_{K_{kh}}(z):z\in M_{kh}-\{0\}\}.
	\end{align*}
	By second variation and Taylor expansion, we have
	\begin{align*}
	P(F;H_R)-P(K_{kh};H_R)=\frac{t^2}{2}\int_{M_{kh}}|\na^{M_{kh}}\vp|^2-|\II_{M_{kh}}|^2\vp^2 d\hi^{m-1}+O(t^3).
	\end{align*}
	Calculating the integral directly by pulling back the volume form on $\R^m$, we have
	\begin{align*}
	\int_{K_{kh}\De F}\dist(z,M_{kh})dz&=(1+O(t))\int_{M_{kh}}^{}d\hi^{m-1}(z)\int_0^{t|\vp(z)|}sds\\&=\frac{t^2}{2}\int_{M_{kh}}\vp^2d\hi^{m-1}+O(t^3).
	\end{align*}
	For details, please refer to Lemma 3.1 in \cite{phi}. Putting these two, and letting $t\to 0,$ we deduce that
	\begin{align}
	\int_{M_{kh}}|\na^{M_{kh}}\vp|^2-|\II_{M_{kh}}|^2\vp^2 d\hi^{m-1}\ge \frac{c_{k,h}}{R^2}\int_{M_{kh}}\vp^2d\hi^{m-1}.
	\end{align}
	To extend (4) to all $\phi\in C^1(M_{kh}),$ let $\psi_j$ be a sequence of cut-off functions so that $\textnormal{spt}\psi_j\s B_{2/j}^m$ and $\psi_j=1$ on $B_{1/j}^m$ with $|D\psi_j|\le C_mj$ everywhere, where $C_m$ is a positive constant depending only on $m$. We know that $\hi^{m-1}(M_{kh}\cap B_r^m)\le c(m)r^{m-1}$ for some constant $c(m)$ depending only on $m$ and $|\II_{M_{kh}}|\le \frac{C}{|z|}$ for some constant $C$ depending only on $k,h.$ Combining these estimates, we can see that the integrand on the left hand side of (4) is dominated by $O(\frac{1}{|z|^2}),$ and thus the integral on the left hand side converges as $j\to\infty$. Let $j\to\infty$ and use dominated convergence. We deduce that (4) is true for all $\vp\in C^1(M_{kh}).$\qed
	\subsection{Proof of Theorem 1}
Define
\begin{align*}
p(z)=&\bigg|\frac{|x|}{\sqrt{k-1}}-\frac{|y|}{\sqrt{h-1}}\bigg|.
\end{align*}	By Lemma 1 and Lemma 2, we have
	\begin{align*}
	|K_{kh}\De F|&\le |(K_{kh}\De F)\cap\{p>\e\}|+|H_R\cap\{p<\e\}|\\
	&\le \int_{(K_{kh}\De F)\cap\{p>\e\}}\frac{p(z)}{\e}\frac{R^2}{|z|^2}dz+|H_R\cap\{p<\e\}|\\
	&=\frac{lR^2}{\e}
	\int_{(K_{kh}\De F)\cap\{p>\e\}}\frac{\dist(z,M_{kh})}{|z|^2}dz+|H_R\cap\{p<\e\}|\\
	&\le \frac{lR^2}{c_{k,h}\e}\int_{(K_{kh}\De F)\cap\{p>\e\}}|\di g|dz+|H_R\cap\{p<\e\}|\\
		&\le \frac{lR^2}{c_{k,h}\e}\bigg(P(F;H_R)-P(K_{kh};H_R)\bigg)+|H_R\cap\{p<\e\}|,
		\end{align*}
	where $l=\sqrt{\frac{1}{h-1}+\frac{1}{k-1}}$ by elementary geometry. Now, we need to get a suitable upper bound for $|H_R\cap\{p<\e\}|$. We have
	\begin{align*}
	|H_R\cap\{p<\e\}|=&\int_{B_R^h}\hi^k\bigg(\bigg\{x\in B_R^k:\frac{|y|}{\sqrt{h-1}}-\e<\frac{|x|}{\sqrt{k-1}}<\frac{|y|}{\sqrt{h-1}}+\e\bigg\}\bigg)dy\\
	&\le \om_k(k-1)^{k/2}\int_{B_R^h}\bigg(\frac{|y|}{\sqrt{h-1}}+\e\bigg)^k-\bigg(\frac{|y|}{\sqrt{h-1}}-\e\bigg)_+^kdy.
	\end{align*}
	We can break down the estimate into two parts, namely
	\begin{align*}
	&\int_{B_{\e\sqrt{h-1}}^h}\bigg(\frac{|y|}{\sqrt{h-1}}+\e\bigg)^k-\bigg(\frac{|y|}{\sqrt{h-1}}-\e\bigg)_+^kdy\\
	=&\int_{B_{\e\sqrt{h-1}}^h}\bigg(\frac{|y|}{\sqrt{h-1}}+\e\bigg)^kdy\\
	\le&2^k\e^{h+k}\om_h(h-1)^{k/2},
	\end{align*}
	and
	\begin{align*}
&	\int_{B_R^h\backslash B_{\e\sqrt{h-1}}^h}\bigg(\frac{|y|}{\sqrt{h-1}}+\e\bigg)^k-\bigg(\frac{|y|}{\sqrt{h-1}}-\e\bigg)^k_+dy\\
=&\int_{B_R^h\backslash B_{\e\sqrt{h-1}}^h}\bigg(\frac{|y|}{\sqrt{h-1}}+\e\bigg)^k-\bigg(\frac{|y|}{\sqrt{h-1}}-\e\bigg)^kdy\\
\le&\frac{1}{(h-1)^{k/2}}\int_{B_R^h\backslash B_{\e\sqrt{h-1}}^h}|y|^k\bigg(\bigg(1+\frac{\e\sqrt{h-1}}{|y|}\bigg)^k-\bigg(1-\frac{\e\sqrt{h-1}}{|y|}\bigg)^k\bigg)dy\\
\le& \frac{2^k}{(h-1)^{k/2}}\int_{B_R^h\backslash B_{\e\sqrt{h-1}}^h}|y|^k\frac{\e\sqrt{h-1}}{|y|}dy\\
\le& \frac{2^k\e}{(h-1)^{(k-1)/2}}\int_{\e\sqrt{h-1}}^{R}r^{k-1}\hi^{m-1}(S_r^{h-1})dr\\
\le&\frac{2^kh\om_h\e}{(h-1)^{(k-1)/2}}\int_{\e\sqrt{h-1}}^Rr^{k+h-2}dr\\
\le&\frac{2^kh\om_h\e}{(h-1)^{(k-1)/2}(m-1)}R^{m-1},
	\end{align*}where we use $(1+t)^k-(1-t)^k\le 2^kt$ for $t\in(0,1),k\in \N.$
Combining the two parts, we have
\begin{align*}
|H_R\cap\{p<\e\}|\le &2^k\om_k\om_h(k-1)^{k/2}(h-1)^{h/2}\e\bigg(\e^{m-1}+\frac{hR^{m-1}}{(h-1)^{(m-1)/2}(m-1)}\bigg)
\end{align*}
Now, note that $\om_j<6$ for all $2\le j\le 11,$ so by substituting the explicit values for $c_{k,h},$ we have
		\begin{align}
		|K_{kh}\De F|\le& \frac{2\times 11^5\sqrt{11}R^2}{\e}\bigg(P(F;H_R)-P(K_{kh};H_R)\bigg)\nonumber\\&+2^{11}6^210^{11/2}2^{3/2}\e(\e^{m-1}+\frac{3}{6}R^{m-1}) \\
		\le& 7\times 10^{10}\bigg(\frac{R^2}{\e}\big(P(F;H_R)-P(K_{kh};H_R)\big)+\e(\e^{m-1}+R^{m-1})\bigg).
		\end{align}
	Let
	\begin{align*}
	 \ai&=\frac{\lb^m(K_{kh}\De F)}{R^m},\\ \de&=\frac{P(F;H_R)-P(K_{kh};H_R)}{R^{m-1}}.
	\end{align*}
	Note that $\ai\le R^{-m}\lb^m(H_R)=\om_k\om_h\le 6^2.$ If $\de\ge 6^2,$
 then $\ai\le\om_k\om_h\le 6\sqrt{\de}.$ Thus we assume $\de \le 6^2.$ Inequality (6) implies
 \begin{align}
 \ai\le 7\times 10^{10}\bigg(\frac{R}{\e}\de+\frac{\e}{R} \big((\e/R)^{m-1})+1\big)\bigg)
 \end{align}
If $\e<\sqrt[13]{35}R,$ then inequality (7) implies
\begin{align*}
\ai\le 7\times 10^{10}(\frac{R}{\e}\de+36\frac{\e}{R} ).
\end{align*}
Note that
\begin{align*}
\frac{R}{\e}\de+36\frac{\e}{R}\ge 12\sqrt{\de}
\end{align*}
	with equality if and only if $\e=R\sqrt{\frac{\de}{36}}.$ Since $\frac{\de}{36}\le1,$ we can let $\e=R\sqrt{\frac{\de}{36}},$ and deduce that
	\begin{align*}
	\ai\le 7\times 12\times 10^{10}\sqrt{\de}.
	\end{align*}\qed

	\section{Proof of Lemma 2}
\subsection{Calculating $\di g$ on $K_{kh}$}	
To make calculations simpler, let $u=(h-1)|x|^2,v=(k-1)|y|^2.$ First, consider the function
\begin{align*}
f(z)=\frac{1}{4}(u-v)u^d.
\end{align*}
We have
\begin{align*}
\pd_{x_i}f=&\frac{h-1}{2}x_i\big((d+1)u^d-dvu^{d-1}\big),\\
\pd_{x_i}\pd_{x_j}f=&\frac{h-1}{2}\de_{ij}\big((d+1)u^d-dvu^{d-1}\big)+(h-1)^2 x_ix_j\big((d+1)du^{d-1}-d(d-1)vu^{d-2}\big),\\
\pd_{y_i}f=&-\frac{k-1}{2}y_iu^d,\\
\pd_{y_j}\pd_{y_i}f=&-\frac{k-1}{2}\de_{ij}u^d,\\
\pd_{y_j}\pd_{x_i}f=&-d(h-1)(k-1)u^{d-1}x_iy_j.
\end{align*}
This gives
\begin{align*}
|\na f|^2=&\frac{h-1}{4}u\big((d+1)u^d-dvu^{d-1}\big)^2+\frac{k-1}{4}vu^{2d},\\
\De f=&\frac{(h-1)k}{2}\big((d+1)u^d-dvu^{d-1}\big)-\frac{(k-1)h}{2}u^d\\&	+(h-1)u\big((d+1)du^{d-1}-d(d-1)vu^{d-2}\big),1\\
(\pd_{x_i}f)(\pd_{x_j}f)(\pd_{x_i}\pd_{x_j}f)=&\frac{(h-1)^2}{8}u\big((d+1)u^d-dvu^{d-1}\big)^3\\&+\frac{(h-1)^2}{4}u^2\big((d+1)u^d-dvu^{d-1}\big)^2\big((d+1)du^{d-1}-d(d-1)vu^{d-2}\big),\\
(\pd_{y_i}f)(\pd_{y_j}f)(\pd_{y_i}\pd_{y_j}f)=&-\frac{(k-1)^2}{8}vu^{3d},
\\
(\pd_{x_i}f)(\pd_{y_j}f)(\pd_{x_i}\pd_{y_j}f)=&\frac{(h-1)(k-1)}{4}du^{2d	}v\big((d+1)u^{d}-dvu^{d-1}\big).\end{align*}
Thus, we have
\begin{align*}
|\na f|^3\di g=&|\na f|^3\di\frac{\na f}{|\na f|}\\=&|\na f|^2\De f-(\pd_{x_i}f)(\pd_{x_j}f)(\pd_{x_i}\pd_{x_j}f)-(\pd_{y_i}f)(\pd_{y_j}f)(\pd_{y_i}\pd_{y_j}f)\\&-2(\pd_{x_i}f)(\pd_{y_j}f)(\pd_{x_i}\pd_{y_j}f)\\
=&\frac{ (h-1) (k-1)}{8} u^{3d-2} (u - v) \bigg((1 + d)^2 \big(-1 + d (-1 + h)\big) u^2\\& + 
d (-2 + d (1 + 2 d - 2 (1 + d) h) + k) u v + d^3 (-1 + h) v^2\bigg).
\end{align*}
\subsection{Calculating $\di g$ on $K_{kh}\cp$}
Now, define
 $$f(z)=\frac{1}{4}(u-v)v^d,$$ which can be obtained by interchanging $u,v$ and $h,k$ and adding an additional minus sign to $f$ in the previous subsection. Thus, by symmetry or by direct computations, we must have
\begin{align*}
|\na f|^2=& \frac{h-1}{4} u v^{2 d} + 
\frac{k-1}{4} v (d u v^{d-1} - (d+1) v^d)^2,\\
|\na f|^3\di g=&|\na f|^3\di\frac{\na f}{|\na f|}\\=&\frac{(h-1)(k-1)}{8} (u - v) v^{3d-2} \bigg(d^3 (k-1) u^2 \\&+ 
d \big(-2 + d + 2 d^2 + h - 2 d (1 + d) k\big) u v\\& + (d+1)^2 (-1 + 
d (k-1)) v^2\bigg).
\end{align*}

Note that if we set $g=\frac{\na f}{|\na f|},$ then $g$ is clearly continuous, and smooth except on $M_{kh}.$ Calculations can show that the derivative of $g$ is of order $O(|z|^{-1})$ near origin, so $g\in W^{1,1}_{\textnormal{loc}}(\R^m,\R^m).$
	
		\subsection{The Cases $(2,k)$}
			We use the basic inequalities $\max\{|x|,|y|\}\le z\le \sqrt{2}\max\{|x|,|y|\}$ and $(a^q+b^q)^{1/q}\le (a^2+b^2)^{1/2}$  for $a,b>0,q\ge 2.$ Also note that $$d(z,M_{kh})=\frac{|\sqrt{u}-\sqrt{v}|}{\sqrt{h+k-2}}$$ by elementary geometry.
	
	If $u>v,$ then choosing $d=3/2,$ we have
	\begin{align*}
	\di  g=\frac{\frac{1}{64} (-1 + k) u^{5/2} (u - v) (25 u^2 + 12 (-11 + k) u v + 27 v^2)}{\bigg(\frac{1}{16} u^2 (25 u^2 + 2 (-17 + 2 k) u v + 9 v^2)\bigg)^{3/2}}.
	\end{align*}
	Let $p_2(t)=27 t^2-48 t+25.$ We have $\min_{[0,1]}p_2=p_2(8/9)=11/3.$
	
	This gives
	\begin{align*}
	\di  g\ge &(k-1)(\sqrt{u}-\sqrt{v})\frac{\sqrt{u}+\sqrt{v}}{\sqrt{u}}\frac{25 u^2 + 12 (-11 + 7) u v + 27 v^2}{ (25 u^2 + 2 (-17 + 2 k) u v + 9 v^2)^{3/2}}\\
	\ge&(k-1)(\sqrt{u}-\sqrt{v})\frac{u^2 p_3(v/u)}{ (25 u^2 + 10 u v + 9 v^2)^{3/2}}\\
	\ge&(k-1)(\sqrt{u}-\sqrt{v})\frac{\frac{11}{3}(|z|/\sqrt{2})^4 }{ (44u^2)^{3/2}}\\
		\ge&(k-1)(\sqrt{u}-\sqrt{v})\frac{\frac{11}{3}(|z|/\sqrt{2})^4 }{ (44 |z|^4)^{3/2}}\\
	\ge&\frac{(k-1)}{2^5 3\sqrt{11}}\frac{\sqrt{u}-\sqrt{v}}{|z|^2}\\
	\ge&\frac{\sqrt{11}}{2^53}\frac{\sqrt{u}-\sqrt{v}}{|z|^2}.	\end{align*}
	
	If $u<v,$ choosing $d=1,$ we have
	\begin{align*}
	\di  g=\frac{\frac{1}{8} (k-1) (u - v) v ((k-1) u^2 + (3 - 4 k) u v + 4 (-2 + k) v^2)}{\bigg(\frac{1}{4} v ((k-1) (u - 2 v)^2 + u v)\bigg)^{3/2}}.
	\end{align*}
	Let $q_2(t)=(k-1)t^2+(3-4k)t+4(k-2).$ We know that $\min_{[0,1]}q_2=q_2(1)=k-6.$
	This gives
	\begin{align*}
	|\di  g|\ge&(k-1)|\sqrt{u}-\sqrt{v}|\frac{\sqrt{u}+\sqrt{v}}{\sqrt{v}}\frac{v^2q_2(u/v)}{((k-1)4v^2+v^2)^{3/2}}\\
	\ge &(k-1)|\sqrt{u}-\sqrt{v}|\frac{(k-6)(k-1)^2(|y|)^4}{(4(k-1)^3+(k-1)^2)^{3/2}|y|^6}\\
	\ge&(k-1)|\sqrt{u}-\sqrt{v}|\frac{(k-6)(k-1)^2(|z|/\sqrt{k})^4}{(4(k-1)^3+(k-1)^2)^{3/2}|z|^6}\\
	\ge&\frac{(k-6)(k-1)^3}{k^2(4(k-1)^3+(k-1)^2)^{3/2}}\frac{|\sqrt{u}-\sqrt{v}|}{|z|^2}\\
	\ge&\frac{6^3}{11^2(4\times 10^3+10^2)^{3/2}}\frac{|\sqrt{u}-\sqrt{v}|}{|z|^2}\\\ge&\frac{1}{11^5}\frac{|\sqrt{u}-\sqrt{v}|}{|z|^2},
	\end{align*}
	where we use $v>u$ if and only if $|x|<\sqrt{k-1}|y|$ and thus $|z|<\sqrt{k}|y|$
	
	This gives $$|\di  g|\ge\frac{1}{11^5} \frac{|\sqrt{u}-\sqrt{v}|}{|z|^2}=\frac{1}{11^5\sqrt{11} }\frac{\dist(z,M_{2,k})}{|z|^2}.$$ 
	
	\subsection{The Case $(h,k)=(3,5)$}

	Choose $d=3/4.$	If $u>v,$ we have
	\begin{align*}
	\di  g=\frac{\frac{1}{32} u^{1/4} (u - v) (49 u^2 - 72 u v + 27 v^2)}{\bigg(\frac{1}{32} \sqrt{u} (49 u^2 - 10 u v + 9 v^2)\bigg)^{3/2}}.
	\end{align*}
Let $$p_3(t)=27 t^2-72 t+49.$$
	We know that $\min_{[0,1]}p_3=p_3(1)=4.$ This yields
	\begin{align*}
	\di  g\ge
	& 4\sqrt{2}(\sqrt{u}-\sqrt{v})\frac{\sqrt{u}+\sqrt{v}}{\sqrt{u}}\frac{u^2p_3(v/u)}{\bigg(49 \times 4 |x|^4+9\times 16 |y|^4\bigg)^{3/2}}\\\ge& 4\sqrt{2}(\sqrt{u}-\sqrt{v})\frac{4 u^2}{\bigg((49\times 4)(|x|^4+|y|^4)\bigg)^{3/2}}\\
	\ge&4\sqrt{2}(\sqrt{u}-\sqrt{v})\frac{4\times 4(|z|/\sqrt{2})^4}{14^3 |z|^6}\\
	\ge&\frac{2\sqrt{2}}{7^3}\frac{\sqrt{u}-\sqrt{v}}{|z|^2}.
	\end{align*}
	
	If $u<v,$ we have
	\begin{align*}
	\di  g=\frac{\frac{1}{16} (u - v) v^{1/4} (27 u^2 - 123 u v + 98 v^2)}{\bigg(\frac{1}{16} \sqrt{v} (9 u^2 - 34 u v + 49 v^2)\bigg)^{3/2}}.
	\end{align*}
	Let $q_3(t)=27 t^2-123 t+98.$ We have $\min_{[0,1]}q_3(t)=q_3(1)=2.$
	This gives
	\begin{align*}
	|\di  g|\ge
	&4|\sqrt{u}-\sqrt{v}|\frac{\sqrt{u}+\sqrt{v}}{\sqrt{v}}\frac{2v^2q_3(u/v)}{\bigg(9\times 4|x|^4+49\times 16|y|^4\bigg)^{3/2}}\\\ge&4|\sqrt{u}-\sqrt{v}|\frac{2v^2}{\bigg(49\times 16(|x|^4+|y|^4)\bigg)^{3/2}}\\
	\ge&4|\sqrt{u}-\sqrt{v}|\frac{2 \times 4^2|y|^4}{28^3|z|^6}\\
	\ge&4|\sqrt{u}-\sqrt{v}|\frac{2 \times 4^2(|z|/\sqrt{3})^4}{28^3|z|^6}\\
	\ge&\frac{2}{3^27^3}\frac{|\sqrt{u}-\sqrt{v}|}{|z|^2},
	\end{align*}where we use $u<v\eq 2|y|^2>|x|$ and thus $|z|^2<3|y|^2.$
	This yields
	\begin{align*}
	|\di  g|\ge \frac{2}{3^2 7^3}\frac{|\sqrt{u}-\sqrt{v}|}{|z|^2}=\frac{\sqrt{3}}{21^3}\frac{\dist(z,M_{5,3})}{|z|^2}.
	\end{align*}

	\section*{Acknowledgements}
	The author is very fortunate to have been introduced to the world of minimal surfaces by Professor Hubert Bray, who has taught the author a lot about geometry and life. The author can't thank him enough for his unwavering support and constant encouragement. Also, he would like to thank Professor Guido De Philippis and Professor Francesco Maggi for many helpful discussions and comments on the problem. He also thanks Professor Robert Bryant, Professor Mark Stern, and Professor Richard Hain for countless helpful discussions. Last but not least, he is indebted to Professor David Kraines, who funded the research in the paper with PRUV
	Fellowship.

\end{document}